\documentclass[11pt,reqno,intlimits]{amsart}
\usepackage[cp1251]{inputenc}
\usepackage[T2A]{fontenc}
\usepackage[russian]{babel}
\usepackage{amsfonts,amsmath,amsxtra,amsthm,amssymb,latexsym}
\usepackage[dvips]{graphicx}
\usepackage{graphicx}
\def\ba{\begin{eqnarray}}
 \def\ea{\end{eqnarray}}
 \def\bt{\begin{theorem}}
 \def\et{\end{theorem}}
\graphicspath{{pictures/}}
\DeclareGraphicsExtensions{.pdf,.png,.jpg}
\usepackage{float}
\usepackage{wrapfig}
\def\be{\begin{equation}}
\def\ee{\end{equation}}
\begin{document}
 
 УДК 51 (091)\qquad \qquad \qquad \qquad \qquad \qquad \qquad \qquad \qquad \qquad \texttt{DOI XX.XXXXX/XXXX-XXXX-XXXX-XXXXX}
 
 \hspace{10mm}
 
 \textit{Посвящается 100-летнему юбилею уравнения Некрасова}
 
 \hspace{10mm}
 
\begin{center}
{\textbf{К ИСТОРИИ УРАВНЕНИЯ НЕКРАСОВА}}\footnote{Основные результаты работы обсуждались на Воронежской весенней математической школе  «Понтрягинские чтения – XXXII» 04.05.2021 [Богатов, 2021].} \\

\textbf{Е. М. Богатов} \\

{\it Статья представлена членом редакционной коллегии  }
\vskip 0.2cm

\end{center}
\begin{center}
Филиал Национального исследовательского технологического университета "МИСиС" \, в г. Губкине Белгородской области, \\
Россия, 309180, Белгородская обл., г. Губкин, Комсомольская улица, 16. \\

Старооскольский технологический  институт им. А.А. Угарова (филиал) Национального исследовательского технологического университета "МИСиС" \, , \\ 
Россия, 309516, Белгородская обл., г. Старый Оскол, мкр. Макаренко, 42. \\ 

\vskip 0.2cm
 {\small E-mail: embogatov@inbox.ru} 
\end{center}

\noindent {\bf Аннотация. }
Появившись в 1921 г. как уравнение волн малой амплитуды на поверхности бесконечно глубокой жидкости, уравнение Некрасова быстро стало источником получения новых результатов. Это проявилось как в области математики (теория нелинейных интегральных уравнений А.И. Некрасова; 1922, позже – Н.Н Назарова; 1941), так и в области механики (переход к жидкости конечной глубины – А.И. Некрасов; 1927 и отказ от малости амплитуды волн – Ю.П. Красовский; 1960). 

Основная задача автора – выяснить предысторию возникновения уравнения Некрасова и проследить изменение подходов к его решению в контексте развития нелинейного функционального анализа  1940-х – 1960-х гг. Пристальное внимание будет уделено вкладу европейских и отечественных математиков и механиков: А.М. Ляпунова, Э. Шмидта, Т. Леви-Чивита, А. Вилля, Л. Лихтенштейна, М.А. Красносельского, Н.Н. Моисеева, В.В. Покорного и др. 

В контексте развития качественных методов исследования уравнения Некрасова будет также освещён вопрос о взаимодействии воронежской школы нелинейного функционального анализа под руководством профессора М.А. Красносельского и ростовской школы нелинейной механики под руководством профессора И.И. Воровича.

\textbf{Ключевые слова:}
{нелинейные интегральные уравнения, уравнение Некрасова, уравнение Гаммерштейна, теория Ляпунова-Шмидта, Т. Леви-Чивита, Н.Н. Назаров, М.А. Красносельский, Ю.П. Красовский, воронежская школа нелинейного функционального анализа; ростовская школа нелинейной механики.}

\smallskip

\noindent {\bf Благодарности:} Работа выполнена при финансовой поддержке РФФИ, проект 20-011-00402. \\
Автор выражает признательность участникам Воронежской весенней математической школы «Понтрягинские чтения – XXXII» за внимание к работе и полезные обсуждения, а также В.П. Богатовой за помощь в доступе к первоисточниками.

\smallskip

\noindent {\bf Для цитирования:} Богатов Е.М. 2021. Об истории уравнения Некрасова.  Прикладная математика \& Физика. 53(4): XX--XX. DOI XX.XXXXX/XXXX-XXXX-XXXX-XXXXX.

\vskip - 0.5cm
\phantom{.}\hskip -0.5cm \underline{\underline{\phantom{AAAAAAAAAAAAAAAAAAAAAAAAAAAAAAAAAAAAAAAAAAAAAAAAAAAAAAAAAAAA}}}
\vskip 0.5 cm

\begin{center}
{
\textbf{ON THE HISTORY OF THE NEKRASOV'S EQUATION}}
\medskip

{\bf Egor Bogatov}
\vskip 0.2cm

{\it Article submitted by a member of the editorial board }
\vskip 0.2cm

{\small 
Branch of National Research University of Science and Technology "MISIS" \, in Gubkin town of Belgorod Region \\
Gubkin, 309180,  Russia;\\
Stary Oskol Technological Institute of National Research University of Science and Technology "MISIS" \, \\
Stary Oskol, 309516, Russia \\

}

{\small  E-mail: embogatov@inbox.ru}

Received Xxxx, xx, 2021
\end{center}
{\small
\noindent {\bf Abstract.}
Appearing in 1921 as an equation for small-amplitude waves on the surface of an infinitely deep liquid, the Nekrasov’s equation quickly became a source of new results. This manifested itself both in the field of mathematics (theory of nonlinear integral equations of A.I. Nekrasov; 1922, later - of N.N. Nazarov; 1941), and in the field of mechanics (transition to a fluid of finite depth - A.I. Nekrasov; 1927 and refusal on the smallness of the wave amplitude - Yu.P. Krasovskii; 1960).

 The main task of the author is to find out the prehistory of the Nekrasov’s equation and to trace the change in approaches to its solution in the context of the nonlinear functional analysis development in the 1940s - 1960s. Close attention will be paid to the contribution of European and Russian mathematicians and mechanics: A.M. Lyapunov, E. Schmidt, T. Levi-Civita, A. Villat, L. Lichtenstein, M.A. Krasnoselskii, N.N. Moiseev, V.V. Pokornyi, etc.

In the context of the development of qualitative methods for the Nekrasov’s equation investigating, the question of the interaction between Voronezh school of nonlinear functional analysis under the guidance of Professor M.A. Krasnoselskii and Rostov school of nonlinear mechanics under the guidance of Professor I.I. Vorovich.

}
\smallskip

\noindent {\bf Acknowledgements:} The work is supported by RFBR, project No 20-011-00402.\\

\noindent {\bf For citation:} Bogatov E.M. 2021. On the history of the Nekrasov’s equation. Applied Mathematics \& Physics. 53(4): XX--XX. (in Russian) DOI XX.XXXXX/XXXX-XXXX-XXXX-XXXXX.

\vskip 0.3cm
\noindent\rule{\textwidth}{0.5pt}
\vskip 0.3cm
{\bf 1. Введение. }
 
 С момента опубликования Александром Ивановичем Некрасовым своего уравнения, описывающего распространение малых волн на поверхности тяжёлой жидкости [Некрасов, 1921]:
 \ba \label {f1}
 \Phi(\theta)=-\frac\mu{12\pi}\int_0^{2\pi}\frac{\sin\Phi(\varepsilon)}{1+\mu\int_0^\varepsilon\sin\Phi(\alpha)d\alpha}K(\varepsilon,\theta)d\varepsilon,
 \ea
 где $K(\varepsilon,\theta)=\ln\left|{\frac{1-\cos(\varepsilon-\theta)}{1-\cos(\varepsilon+\theta)}}\right|$; $\mu$ - параметр, $\Phi(\theta)$ - угол наклона касательной к профилю волны, прошло 100 лет. Несмотря на это, полноценный историко-научный анализ его математических работ, в особенности по нелинейным интегральным уравнениям, до 2020 г. проведён не был (некоторые фрагменты такого анализа можно найти в [Вайнберг, Треногин, 1962, §2], [Секерж-Зенькович, 1960], [Волгина, Тюлина, 2001, Гл. 3]). В 2021 вышла статья российского учёного Н.Г. Кузнецова «Повесть о двух интегральных уравнениях Некрасова» [Kuznetsov, 2021], в которой указанный пробел был частично восполнен. В настоящей работе мы попытаемся дополнить исследование Кузнецова новыми фактами и предложить свой взгляд на историю уравнения (\ref{f1}).

Основным инструментом для построения решений (\ref{f1}), предложенным Некрасовым, было разложение в степенной ряд. С того времени появились новые подходы к исследованию подобных уравнений, основанные, прежде всего, на методах нелинейного функционального анализа. При этом само уравнение не потеряло своей актуальности как для математиков, так и для механиков [Kuznetsov, 2021, § 3].
 
Задача автора – проследить эволюцию подходов к изучению уравнения (\ref{f1}), ответив на следующие вопросы:
\begin{enumerate}
	\item  На каком уровне находилась теория нелинейных интегральных уравнений до Некрасова?
	\item  Что привело Некрасова к нелинейным интегральным уравнениям? На какие работы он опирался? 
	\item  Проводились ли подобные исследования в других странах?
	\item  Можно ли сопоставить метод Некрасова с другими методами?
	\item  Кто и когда занимался продолжением идей Некрасова до начала 1960-х гг.?
\end{enumerate}

\smallskip

{{\bf 2. Предыстория : работы А.М. Ляпунова и Э. Шмидта.}

Прежде, чем перейти к исследованиям самого Некрасова, обозначим уровень развития нелинейных интегральных уравнений к концу 1910-х гг. Одним из первых таких уравнений было уравнение (\ref{f2}), возникшее при отыскании форм равновесия вращающихся жидкостей в работе 1903 г. выдающегося русского математика Александра Михайловича Ляпунова [Liapunoff, 1903]
\ba \label {f2}
\int_0^A\rho^2 (a^\prime)da^\prime\int_S(1+(\varsigma^\prime))^2\left(1+\frac{\partial(a^\prime\varsigma^\prime)}{\partial a^\prime}\right)\frac{d\sigma^\prime}{\Delta}+\frac{\omega^3}{2g}a^2(1+\varsigma^2)\sin^2\theta=F(a).
\ea
Здесь $\omega$ – угловая скорость вращения жидкой массы,  g – постоянная всемирного тяготения; $\Delta$ - расстояние между переменной точкой $P'(r', \theta' ,\varphi')$  внутри жидкости и фиксированной точкой $P(r, \theta,\varphi)$ на $S$; $(r,\theta,\varphi)$ -  сферические координаты, $\zeta$- искомая функция параметра а  и координат $(\theta,\varphi)$, характеризующая отклонение поверхности уровня жидкой массы от сферы радиуса $а$; $F(a)$ -  заданная функция; интегрирование во втором интеграле производится по поверхности единичной сферы $S$.

Более общий случай был рассмотрен Ляпуновым в работе 1905 г. [Liapounoff, 1905]. В процессе доказательства существования неэллипсоидальных форм равновесия жидкой массы он вывел уравнение
\ba \label {f3}
R\zeta-\int_S\frac{\varsigma'd\sigma}{\Delta}=W(\zeta),
\ea
где $R=R(\theta,\varphi)$ - заданная функция. Основную сложность этого уравнения составляла правая часть, состоящая из бесконечного числа интегро-степенных слагаемых (пояснение будет дано ниже, при описании вклада Э. Шмидта).

Ляпунов искал малые решения (\ref{f3}) в виде ряда по степеням некоторого малого параметра (равного отношению центробежной силы к силе тяжести на экваторе жидкой планеты). Сходимость этого ряда доказывалась методом мажорант. Доказательство наличия нескольких решений опиралось на нетривиальную разрешимость соответствующего линейного интегрального уравнения.

Ляпунов получил свои результаты ещё до того, как теория линейных интегральных уравнений обрела более или менее сложившуюся форму. Его аргументация отличалась большой сложностью, тем не менее его теория малых решений нелинейных интегральных уравнений вызвала в Европе живой интерес. Математический аппарат теории уравнений вида (\ref{f3}) был разработан немецким математиком Эрхардом Шмидтом (1907-1908), учеником Д. Гильберта [Schmidt, 1907] -- [Schmidt, 1908]. Он, как и Ляпунов, ограничился случаем малых решений уравнений
\ba \label {f4}
u(x)-\int_a^bK(x,s)u(s)ds=f(x)+\int_a^bK_1(x,s)v(s)ds-\sum_{m+n\ge2}U_{mn}\begin{pmatrix} x \\ u, v \end{pmatrix},
\ea
где $x\in(a,b)$; $v(s)$, $f(x)$, $K(x,s)$ и $K_1(x,s)$ - заданные функции; $u(x)$ - искомая функция;
$U_{mn}\begin{pmatrix} x \\ u, v \end{pmatrix}$ - конечная или бесконечная сумма интегро-степенных слагаемых. При малых $|u|$, $|v|$ ряды правой части (\ref{f4}) должны, по условию, сходиться абсолютно и равномерно при всех $x\in(a,b)$. \\
В качестве простейшего примера уравнения (\ref{f4}) можно предложить уравнение
$$
u(x)-\int_0^1K(x,s)u(s)ds=v(x)+\int_0^1\int_0^1K(x,y_1,y_2)u^2(y_1)v^3(y_2)dy_1dy_2.
$$

Для исследования разрешимости уравнения (\ref{f4}) в малом Шмидт представил его как уравнение с известной правой частью и применил к нему теорию Фредгольма о представлении решений линейного уравнения посредством резольвенты (разрешающего ядра) $R(x,t)$. В соответствии с этой теорией здесь необходимо было выделить два варианта дальнейшего рассмотрения, в зависимости от того, является ли оператор, стоящий в левой части (\ref{f4}) обратимым, или нет.

В случае, когда единица является собственным значением интегрального оператора \textit{А}, решение линейного интегрального уравнения с «замороженной» правой частью будет зависеть от нескольких произвольных постоянных $\xi_i$ (по величине кратности единицы, как собственного значения оператора А). При этом оказалось естественным, ориентируясь на принцип Коха обращения ряда, искать решение исходного уравнения в виде \footnote{В этой формуле $V_n\begin{pmatrix} x \\ v \end{pmatrix}$ - подлежащая определению интегро-степенная форма степени $n$.}.
\ba \label {f5}
u(x)=\sum_{m+n\ge1}\xi^mV_n^m\begin{pmatrix} x \\ v \end{pmatrix}.
\ea

Для определения возможных значений $\xi$ оставалось подставить ряд (\ref{f5}) в выражение для $u$ через резольвенту линейного уравнения (по Фредгольму) и получить уравнение разветвления \footnote{Такой вид уравнения (\ref{f6}) соответствует случаю, когда кратность единицы, как собственного значения оператора A, равна одному.}:
\ba \label {f6}
\xi=\sum_{m+n\ge1}\xi^m\int_a^b\varphi(s)V_n^m\begin{pmatrix} s \\ v \end{pmatrix}ds;
\ea
$\varphi(x)$ - собственная функция линейного оператора А.

Шмидт показал, что в правой части уравнения (\ref{f6}) фактически отсутствуют первые степени $\xi$, значит, в общем случае оно имеет не единственное решение. Независимо от величины модуля функции $v(x)$, число решений исходного уравнения определяется числом решений уравнения разветвления (\ref{f6}). Из уравнения (\ref{f6}) можно также получить информацию о виде всех нетривиальных малых решений уравнения (\ref{f4}).

После опубликования работ Ляпунова и Шмидта стало очевидным, что неединственность решений нелинейных интегральных уравнений является их типичным свойством. Уравнение разветвления быстро превратилось в рабочий инструмент для исследования интегральных уравнений [Богатов, Мухин, 2021].

\smallskip

{\bf 3. Задача об обтекания препятствия - исток уравнения Некрасова.}

 А.И. Некрасов был учеником Н.Н. Жуковского – основателя московской школы механики, имеющей международное признание. Одной из первых работ Некрасова, связанных с уравнением (\ref{f1}), была статья \textit{«О прерывном течении жидкости в двух измерениях вокруг препятствия в форме дуги круга»}  (1922) [Некрасов, 1922], в которой он обобщил результаты итальянского математика и механика Тулио Леви-Чивита\footnote{В 1934 г. Леви-Чивита был избран почётным членом Академии наук СССР.} и его французского коллеги Анри Вилля\footnote{Действительный член Парижской Академии наук с 1932 г.}. Обзор более ранних результатов  по теории струйных течений конца XIX в. – начала XX в. (работы Г. Гельмгольца, Г. Кирхгофа, Н. И Жуковского и др.) был сделан учеником А.И. Некрасова, механиком Л.И. Седовым [Седов, 1939, § 10].

К Леви-Чивита восходит исследование подобных задач на плоскости методами комплексного анализа\footnote{Метод Леви-Чивита описан в книге М.И. Гуревича [Гуревич, 1979, Гл. IV, § 17].}. При его подходе предполагалось, что решение будет зависеть от нахождения некоторой функции $\varphi$, способ определения которой не был  указан [Levi-Civita, 1907]. Кроме того, данная функция никак не была связана с формой препятствия текущей жидкости. Вилля взялся исправить ситуацию, продолжив исследования Леви-Чивиты по совету своего учителя, французского механика М. Бриллюэна [Tazzioli, 2017]. Вилля применил конформное отображение области с неизвестной границей на полукольцо и ввёл новую функцию $\theta$, которая уже имела тесную связь с формой препятствия. В своей диссертации 1911 г., сведя задачу к отысканию двух сопряжённых тригонометрических рядов переменной $\theta$, он получил некоторое соотношение между ними, которое привело его к нелинейному интегро-дифференциальному уравнению вида\footnote{Вывод уравнения (\ref{f7}) можно найти в статье [Гуревич, 1970, § 2].} [Villat, 1911]:
\ba \label {f7}
\theta'(\sigma)=\lambda H(\theta,\sigma)\exp\left\{-\frac1\pi\int_0^\pi\theta'(\varepsilon)K(\varepsilon,\sigma)d\varepsilon\right\},
\ea
где $\lambda>0$ - некоторая постоянная; $H(\theta,\sigma)$ - некоторая заданная функция; $K(\varepsilon,\sigma)=\ln\left|\frac{\sin\frac{\varepsilon-\sigma}{2}}{\sin\frac{\varepsilon+\sigma}{2}}\right|$.

Отметим, что до начала 1920-х гг. методы решения уравнения (\ref{f7}) отсутствовали. Они были предложены в работах А.И. Некрасова применительно к аналогу уравнения (\ref{f7}) и основывались на использовании функции Жуковского (логарифма от комплексной скорости течения жидкости $\frac{dw}{dz}$; $w=\varphi + i\psi$ – комплексный потенциал). Задачу обтекания криволинейной дуги Некрасов свёл к уравнению относительно некоторой вспомогательной функции $h(\theta)$, зная которую, можно найти функцию Жуковского. Следуя методу Вилля, он привёл задачу к отысканию двух сопряжённых тригонометрических рядов
\ba \label {f9}
g(\theta)=\sum_{n=1}^\infty a_{2n+1}\cos(2n+1)\theta;
\ea
\ba \label {f10}
h(\theta)=\sum_{n=1}^\infty a_{2n+1}\sin(2n+1)\theta,
\ea
между которыми имелось следующее соотношение
\ba \label {f11}
g'(\theta)=-\frac{\pi\lambda}{2}e^{-h(\theta)}\sin\theta(1-\sin\theta).
\ea
Здесь $\lambda>0$ -- некоторая постоянная, $\theta$ – аргумент параметрического переменного $u$, от которого зависят $w$ и $\frac{dw}{dz}$ (определение этой зависимости позволяет рассчитать все силовые, кинематические и геометрические характеристики течения).

Соотношение (\ref{f11}) привело Некрасова к нелинейному интегральному уравнению
\ba \label {f12}
h(\theta)=F(\theta)+\lambda\int_0^\frac{\pi}{2}\left[1-e^{-h(\varepsilon)}\right]\frac{(1+2\sin\varepsilon-\cos2\varepsilon)}{4}\ln\left|\frac{\tg\frac{\varepsilon-\theta}{2}}{\tg\frac{\varepsilon+\theta}{2}}\right|d\varepsilon,
\ea
где $F(\theta)$ – некоторая известная функция.

Решение уравнения (\ref{f12}) Некрасов представил рядом по степеням $\lambda$, сходящимся при $\lambda<\lambda_0$ (при этом числу $\lambda_0$ соответствует дуга окружности, ограничивающей препятствие, равная $20^{\circ}$).

\smallskip

{\bf 4. Задачи плоской теории волн.}

Задачам нелинейной теории волн установившегося вида на поверхности тяжёлой жидкости Некрасов посвятил семь работ, вышедших в 1919 -- 1928 гг.[Некрасов, 1919] -- [Некрасов, 1928]. На протяжении XX в. подобные задачи решались приближённо, из-за сложности этих задач в них отбрасывались нелинейные члены, что оставляло открытым само существование установившихся волн. В работе [Некрасов, 1922a] Некрасов впервые представил алгоритм решения указанной задачи для волн на поверхности бесконечно глубокой  тяжёлой жидкости. Как и в приведённой выше задаче об обтекании препятствия, решение задачи о волнах было основано на идее, идущей от Вилля: область, занятая одной волной конформно отображается внутрь единичного круга так, что свободная поверхность жидкости соответствует границе этого круга. При этом опять возникают тригонометрические ряды вида (\ref{f9})-(\ref{f10}) и соотношение между ними вида (\ref{f11}), что в конечном итоге \footnote{Некрасов применил здесь преобразования, идею которых предложил один из основателей Московской математической школы профессор Н.Н. Лузин [Некрасов, 1961, с.44].} даёт уравнение (\ref{f1}) [Некрасов, 1961, с. 255]. 

Для волн малой амплитуды ($sin\Phi\approx\Phi$) уравнение (\ref{f1}) переходит в линейное интегральное уравнение вида
\ba \label {f13}
\Phi(\theta)=-\frac{\mu}{12\pi}\int_0^{2\pi}\Phi(\varepsilon)\ln\left|\frac{1-\cos(\varepsilon-\theta)}{1-\cos(\varepsilon+\theta)}\right|d\varepsilon.
\ea
Некрасов показал, что если параметр $\mu$ принадлежит промежутку $[0, \mu_1^*)$, где $\mu_1^*$ -- первое собственное значение уравнения (\ref{f13}), то установившиеся волны отсутствуют.  Математический смысл этого результата состоит в том, что первое собственное значение линеаризованного интегрального уравнения является бифуркационным значением исходного (нелинейного) интегрального уравнения.

Для исследования уравнения (\ref{f1}), (\ref{f11})  Некрасов разрабатывает самостоятельную теорию [Некрасов, 1922b] -- [Некрасов, 1922c], в которой рассматривает нелинейное интегральное уравнение общего вида \footnote{Специальный случай уравнения, введённого  через 8 лет после Некрасова в работе [Hammerstein, 1930].} с вырожденным ядром:
\ba \label {f14}
f(x)=\lambda\int_a^bf(y)K(x,y)dy+\varepsilon\lambda\int_a^bR[\lambda,y,f(y)]K(x,y)dy.
\ea
Здесь $K(x,y)=\sum_{n=1}^\infty\frac{\varphi_n(x)\varphi_n(y)}{\lambda_n}$, $\{\varphi_i(x)\}$ - некоторая  система ортонормированных на отрезке $[a,b]$ функций; $\lambda_i$  - собственные значения однородного уравнения Фредгольма, получающегося из (\ref{f14}) при подстановке $\varepsilon=0$. 

Некрасов искал решения (\ref{f14}) в виде
\ba \label {f15}
f(x)=\sum_{i=0}^\infty(\lambda')^{p+i}\cdot f_i(x),
\ea
где $\lambda=\lambda_1\pm\varepsilon\cdot\lambda'$ - малый параметр, $p=\frac{1}{m-1}$; $m\geq2$ - наименьшая степень разложения функции $R$ по степеням $f(y)$:
$$
(\lambda_1\pm\varepsilon\cdot\lambda')R[\lambda, y, f(y)]=(\lambda')^{p+1}R_1+(\lambda')^{p+2}R_2+\cdots,\, R_i=R_i(y, \lambda_1, \varepsilon).
$$
Используя специальный приём, Некрасов пришёл к выражению $f_m (x)$ через функции системы ${\varphi_m (x)}$:
$$
f_0(x)=C_0\varphi_1(x),
$$
\ba \label {f16}
f_m(x)=C_m\varphi_1(a)+\sum_{i=2}^\infty\frac{\lambda_1}{\lambda_i-\lambda_1}a_{im}\varphi_i(x), m=1,2,\ldots
\ea
для которого величины $C_m$ и $a_{im}$ определяются однозначно, а вот величина $C_0$ - нет:
$$
C_0=\mp\int_a^bR_1[\lambda_1, y, C_0\varphi(y)]dy
$$
[Некрасов, 1961, с. 87], что приводит к ветвлению.

Для доказательства сходимости рядов (\ref{f15}) - (\ref{f16}) Некрасов  применял (как Ляпунов и Шмидт) метод мажорант. Однако, в отличие от своих предшественников Некрасов решал исходное нелинейное интегральное уравнение \textit{непосредственно} (не используя редукцию задачи к решению преобразованного функционального уравнения).

Для построения решения уравнения (\ref{f1}) Некрасову удалось использовать более простое соображение. Полагая $\mu=3+\lambda$, он отыскивал решение этого уравнения в виде степенного ряда
\ba \label {f17}
\Phi(\theta,\lambda)=\sum_{k=1}^\infty\lambda^k\Phi_k(\theta).
\ea
Подстановка данного ряда в уравнение (\ref{f1}) привела к рекуррентной системе для определения $\Phi_k$:
\begin{equation}\label {f18}
\left\{\begin{gathered}
			\Phi_1(\theta) = 3\int_0^{2\pi}\Phi_1(\varphi)K(\varphi,\theta)d\varphi
			\\
			\Phi_2(\theta) = 3A(\Phi_2)+A_2(\Phi_1)
			\\
			\Phi_3(\theta) = 3A(\Phi_3)+A_3(\Phi_1,\Phi_2)
			\\
			\cdots
			\\
			\Phi_n(\theta) = 3A(\Phi_n)+A_n(\Phi_1,\Phi_2,\ldots,\Phi_{n-1})
			\\
			\cdots
	\end{gathered}\right.
\end{equation}
где $A_k$ – нелинейные операторы по всем своим переменным, кроме последней, но при $k=3,4,\ldots$ операторы $A_k$ являются линейными относительно $\Phi_{k-1}$ [Вайнберг, Треногин, 1962, § 2]. Оператор $A$ – линейный интегральный оператор с ядром $K(\varphi,\theta)$. 

Следует отметить, что за работу [Некрасов, 1921] по теории волн Некрасов был удостоен премии Жуковcкого [Секерж-Зенькович, 1960, с. 154].

В более поздней своей работе Некрасов исследует задачу о волне установившегося вида на поверхности конечной глубины [Некрасов, 1928]. В этом случае область, занятая одной волной, конформно отображается в круговое кольцо. Опираясь опять на метод Вилля, Некрасов сводит задачу к решению нелинейного интегрального уравнения, аналогичного (\ref{f1}), в котором 
$$
K(\varepsilon,\theta)=-\ln\left|\frac{\sigma\left[\frac\omega\pi(\theta+\varepsilon)\right]\sigma_3\left[\frac\omega\pi(\theta-\varepsilon)\right]}{\sigma_3\left[\frac\omega\pi(\theta+\varepsilon)\right]\sigma\left[\frac\omega\pi(\theta-\varepsilon)\right]}\right|
$$

Здесь $\omega=Const$; $\sigma_3$ и $\sigma$ – это эллиптические сигма-функции Вейерштрасса [Некрасов, 1947].

Поскольку основные результаты по теории нелинейных  интегральных уравнений Некрасова были опубликованы в малоизвестных журналах, они не получили широкого  распространения в то время ни в России, ни за рубежом.

 Отметим, что на международном съезде по теоретической механике (1924) доклад Некрасова (по причине отсутствия докладчика) был прочитан Леви-Чивитой [Лапко, Люстерник, 1957, с. 119]. После этого некоторые западноевропейские учёные обратили внимание на статьи Некрасова 1921-22 гг. (см., например, книгу Л. Лихтенштейна [Lichtenstein, 1931, с. 54]).  Леви-Чивита работал над определением профиля устоявшихся безвихревых волн на поверхности жидкости бесконечной глубины параллельно с Некрасовым. Об этом говорит тот факт, что в 1921 г. он включил эту тему в цикл лекций по вопросам классической и релятивистской механики [Levi-Civita, 1922]. В 1925 г. Леви-Чивита опубликовал \textit{своё} «уравнение волны» [Levi-Civita, 1925], [Stoker, 1957, с. 523-524]:
 \ba \label {f19}
 \frac{d\Phi}{d\theta}=pe^{-3\Phi}\sin\Psi,
 \ea
 где $p>0$ - параметр;\\
 $\Psi=Im \, \omega$, $\Phi=Re \, \omega$; $\varsigma=\rho e^{i\theta}$; \\
 $\omega=\omega(\varsigma)$ – аналог потенциала скорости течения жидкости.
  
 Функция $\omega$ предполагается аналитической внутри круга $|\varsigma|<1$, непрерывной при $|\varsigma|\leq1$ и обладающей свойством $\omega(0)=0$. Уравнение (\ref{f19}) относится к точкам окружности $|\varsigma|=1$; оно является следствием условия постоянства давления жидкости на её поверхности.
 
 Существование решения уравнения (\ref{f19}) доказывалось Леви-Чивитой с помощью довольно сложной процедуры, в которой были задействованы нелинейные тригонометрические выражения и метод мажорант. Он оказался в более выгодном, по сравнению с Некрасовым, положении: 50-страничная работа [Levi-Civita, 1925] вышла в престижном немецком журнале \textit{Mathematische Annalen} и быстро стала классической. Но и исследования Некрасова долгое время продолжали оставались актуальными, и в 1951 году на эту тему им была опубликована монография [Некрасов, 1951], за которую Некрасов получил Сталинскую премию. После этого его результаты стали более востребованы западно-европейскими учёными (см., например, [Stoker, 1957, с. IX]) и встал вопрос о приоритете Некрасова в нелинейной теории волн. Для решения этого вопроса в начале 1960-х гг. во Франции было осуществлено подробное сопоставление результатов Некрасова и Леви-Чивита, которое показало их равносильность [Jolas, 1962].
 
 В начале 1940-х гг. метод Некрасова был «переоткрыт\footnote{Назаров не был знаком с результатами Некрасова, поскольку он не ссылается на работы  [Некрасов, 1919] -- [Некрасов, 1928].}» среднеазиатским математиком Николаем Николаевичем Назаровым [Назаров, 1941], учеником В.А. Стеклова. Назаров стал применять указанный метод для абстрактных аналогов уравнений (\ref{f1}), (\ref{f12}):
 \ba \label {f20}
 u(x)=\lambda\int_a^bH(x,y,u(y))dy
 \ea
 с аналитическими операторами $H$, раскладывающихся в ряд в окрестности известного решения $u_0(x)$:
 $$
 H(x,y,u_0(x)+v)=\sum_{k=0}^\infty H_k(x,y)v^k.
 $$
 В этих условиях Назаров преобразовал исходное уравнение к равенству
 $$ 
 v(x)-\lambda_0\int_a^bH_1(x,y)v(y))dy=\int_a^b\left[\mu(H_0+H_1v)+(\mu+\lambda_0)(H_2v^2+H_3v^3+\cdots)\right]dy
 $$
 и представил его малые решения в виде рядов по целым или дробным степеням параметра $\mu=\lambda-\lambda_0$. Коэффициенты этих рядов были определены рекуррентно из бесконечной системы уравнений, а их сходимость установлена путём построения мажорант.
 
 Важным моментом в результатах Назарова является возможность судить о существовании решения уравнения (\ref{f20}) для любого значения $\lambda$, лежащего  в окрестности $\lambda_0$, на основе знания решения $u_0$ исходного уравнения, соответствующего $\lambda_0$. Назаровым были также найдены условия единственности и неединственности такого решения [Назаров, 1945].
 
 Впоследствии под \textit{методом Некрасова-Назарова} стали подразумевать представление решения уравнения в виде рядов по степеням $\lambda-\lambda_0$ или $(\lambda-\lambda_0)^\frac1n$, причём функциональные коэффициенты этих рядов находятся последовательным решением линейных интегральных уравнений, а сходимость полученных рядов устанавливается методом мажорант [Ахмедов, 1957, с.137]. 
 
\smallskip

{\bf 5. Метод Ляпунова-Шмидта в решении уравнения Некрасова.}

 В 1930-м г. немецкий математик А. Гаммерштейн ввёл в рассмотрение уравнение\footnote{Оно получило название \textit{уравнение Гаммерштейна}.}
 \ba \label {f21}
 y(x)=\int_a^bK(x,s)f(s,y(x))ds,
 \ea
 и доказал его разрешимость «в целом», используя вариационный подход [Hammerstein, 1930], [Bogatov, 2020]. Кроме того, он применил к изучению уравнения (\ref{f21}) метод Ляпунова-Шмидта и вывел для него систему уравнений разветвления [Смирнов, 1936, § 15].
 
 Подход Гаммерштейна был взят на вооружение его коллегой Л. Лихтенштейном, который исследовал задачи теории волн в постановке Леви-Чивиты \mbox{[Lichtenstein, 1931, Гл. II, § 1]}. Лихтенштейн ввёл условия нечётности на функции $\Psi$ и чётности функций $\Phi$ из уравнения (\ref{f19}), связанные с предполагаемой симметричностью волн. Это позволило ему, опираясь на теорию потенциала, свести задачу к системе нелинейных интегральных уравнений Вольтерра-Фредгольма
 
 \begin{equation}
 	\left\{\begin{gathered}
 		\Psi(\varphi) = \frac{p}{\pi}\int_0^{2\pi}\ln\frac{1}{\theta}e^{-3\Phi}\sin\Psi d\theta
 		\\
 		\Phi(\varphi) = p\int_0^\varphi e^{-3\Phi}\sin\Psi d\theta
 		\\
 	\end{gathered}\right.
 \end{equation}
 
 и вывести для неё уравнение разветвления [Lichtenstein, 1931, с. 48-52].
 
 Дальнейшее продвижение в исследовании уравнения Некрасова (\ref{f1}) были связаны с развитием методов Ляпунова-Шмидта в контексте общей теории операторов в банаховых пространствах [Вайнберг, Айзенгендлер, 1966, § 5]. 
 
 Эта теория позволяла записать уравнение (\ref{f1}) (и ему подобные) в операторной форме:
 \ba \label {f22}
 \varphi=\lambda A(\varphi,\lambda),
 \ea
 где $\lambda=-\frac{\mu}{12\pi}$; $A$ -- вполне непрерывный оператор:
 $$
 A(\varphi,\lambda)=\int_0^{2\pi}\frac{\sin\varphi(\varepsilon)}{1-12\pi\lambda\int_0^\varepsilon\sin\varphi(\alpha)d\alpha}K(\varepsilon,\varphi)d\varepsilon
 $$
 Для исследования уравнения (\ref{f22}) на предмет наличия малых решений удобно преобразовать его к виду
 \ba \label {f23}
 x=\lambda B(x),
 \ea
 где $B$  – оператор Ляпунова-Шмидта (данное преобразование было выполнено Н.Н. Моисеевым в [Моисеев, 1957]). Тогда теорема существования нелинейных волн не требует специальных доказательств – она является прямым следствием теории Ляпунова-Шмидта. 
 
 Редукция к интегральным уравнениям Шмидта оказалась выгодной и с прикладной точки зрения. Она позволила применить метод Ляпунова-Шмидта для эффективного расчёта нелинейных волн. Как показал Моисеев, уравнение разветвления, соответствующее интегральному уравнению (\ref{f23}), имеет \textit{«полезный»} вид
 $$
 c=f(\lambda,\alpha)
 $$
 $c$ - скорость движения жидкости; $\lambda$ – длина волны; $\alpha$ – амплитуда.
 
 В начале 1960-х метод Ляпунова-Шмидта был также применён к исследованию уравнения (\ref{f22}) непосредственно [Hyers, 1964, c. 322-324].
 
 В середине 1950-х гг. М.А. Красносельский высказал гипотезу о том, что для нелинейных интегральных уравнений (\ref{f20}) формальные решения в виде рядов всегда являются истинными решениями при достаточно малых значениях параметров [Покорный, 1958, с. 711]. Справедливость этой гипотезы была проверена В.В. Покорным, который обобщил теорию формальных рядов с числовыми коэффициентами, развитую С. Бохнером и У. Мартиным, на случай формальных степенных рядов с функциональными коэффициентами [Покорный, 1958]. Рассматривая аналог «обобщённого уравнения Некрасова» (\ref{f14}) в виде
 \ba \label {f24}
 \varphi(x)=\int_0^1A_{10}(x,y)\varphi(y)dy+\int_0^1\Gamma[x,y,\varphi(y);\lambda]dy,
 \ea
 где $\Gamma[x,y,z;\lambda]=\lambda A_{01}(x,y)+\sum_{k+i\geq 2}A_{ki}(x,y)z^k\lambda^i,$
 
 Покорный представил его решение, в соответствии с методом Ляпунова-Шмидта, следующим образом [Покорный, 1956, с. 43]
 \ba \label {f25}
 \varphi(x,\lambda)=\psi(x,\lambda)+\alpha(\lambda)\omega(x).
 \ea
 Здесь $\omega(x)$ – нормированная собственная функция ядра $A_{10}(x,y)$; $\alpha(\lambda)$ определяется условием ортогональности функций $\omega(x)$ и $\psi(x,\lambda)$ в пространстве $L^2(0,1)$.
 
 Функция $\psi(x,\lambda)$ определялась решением уравнения
 \ba \label {f26}
 \psi(x)=\int_0^1 A_{10}(x,y)\psi(y)dy+f[x,\psi(x);\alpha,\lambda],
 \ea
 где через $f[x,\psi(x);\alpha, \lambda]$ обозначен интеграл $\int_0^1\Gamma[x,y,\varphi(y);\lambda]dy$ при подстановке в него выражения для $\varphi(x,\lambda)$ из (\ref{f25}).
 
 Условие разрешимости уравнения (\ref{f26}) даёт зависимость между $\alpha$ и $\lambda$, которую можно представить в упрощённом виде
  \ba \label {f27}
  F(\alpha,\lambda)=0,
  \ea
  где $F(\alpha,\lambda)$ – аналитическая  в некоторой окрестности точки (0,0) относительно обеих своих переменных функция.
  
  Покорный показал, что система уравнений (\ref{f18}), используемая в методе Некрасова-Назаро\-ва, эквивалентна уравнению разветвления, составленному для уравнения (\ref{f24}) и даже совпадает с ним, если правую часть уравнения разветвления разложить в ряд по соответствующему параметру и условие (\ref{f27}) заменить условием обращения в нуль всех коэффициентов этого разложения [Покорный, 1960]. Попутно Покорный обнаружил, что применение методов Ляпунова-Шмидта и Некрасова-Назаро\-ва требуют преодоления одинаковых технических трудностей – построения резольвенты и решения уравнения разветвления.
  
\smallskip

{\bf 6. Качественный анализ уравнения Некрасова.}
  
  Применение качественных методов к анализу задач гидродинамики восходят к французскому учёному Жану Лере (1935) [Leray, 1935a]-[Leray, 1935b], который исследовал вопросы существования и единственности решений задачи о струйном обтекании препятствия неограниченным потоком (см. также [Гуревич, 1979, Гл. IV, § 18]). Лере опирался на созданную им вместе с польским математиком Юлиушем Шаудером теорию степени отображения [Leray, Schauder, 1934] (историю вопроса см., например, в [Mawhin, 2006]). В середине 1950-х гг. подход Лере был распространён на задачи плоской поверхностных теории волн тяжёлой жидкости в постановке Леви-Чивиты (Р. Жербе, [Gerber, 1955]), однако и здесь дальше вопросов разрешимости дело не пошло.
  
  Поворотным моментом следует, по-видимому, считать защиту докторской диссертации советского математика Марка Александровича Красносельского [Красносельский, 1950], в которой идеи Лере-Шаудера получили своё дальнейшее развитие. Красносельский смотрел на нелинейный функциональный анализ, как на возможность проводить \textit{качественное  исследование} уравнений, не прибегая к построению их решений по примеру того, как это происходило в теории дифференциальных уравнений  [Немыцкий, Степанов, 1947]. Он разработал ряд топологических методов анализа, которые нашли своё применение в различных прикладных областях [Красносельский, 2000], в том числе, при анализе задач теории волн.
  
  В работе [Красносельский, 1956] Красносельский рассмотрел оператор Некрасова вида
  \ba \label {f28}
  A_1(\varphi,\mu)=\mu\int_0^{2\pi}\frac{K(x,y)\sin\varphi(y)dy}{1+\mu\int_0^y\sin\varphi(t)dt},
  \ea
  где
  \ba \label {f29}
  K(x,y)=\sum^\infty_{n=1}\frac{\sin nx\sin ny}{\mu_n},
  \ea
  а $\mu_n=\mu_n^\delta$ – положительные собственные значения ядра $K(x,y)$, имеющие различный вид для жидкости конечной $(\delta=1)$ и бесконечной $(\delta=2)$ глубины. Точные формулы для $\mu_n^\delta$ имеют следующий вид: $\mu_n^1={3ncoth\frac{2\pi nh}{\lambda}}$, где $\lambda$ –  длина волны, $h$ – глубина жидкости [Некрасов, 1961, c. 405]; $\mu_n^2={3n}$.
  
  Красносельский показал, что оператор $A_1$ допускает разложение по формуле Тейлора в окрестности нуля $\theta$ банахова пространства $E$:
  \ba \label {f30}
  A_1(\varphi,\mu)=\mu B\varphi + C(\varphi,\mu)+D(\varphi,\mu)
  \ea
  где $B$ – линейный вполне непрерывный оператор (производная Фреше оператора $A_1$ в точке $\theta$);
  
  $C(\varphi,\mu)$ – оператор $k$-го порядка относительно $\varphi$:
  \ba \label {f31}
  C(\alpha\varphi,\mu)=\alpha^kC(\varphi,\mu);
  \ea
  $D(\varphi,\mu)$ - оператор высшего, чем $k$ порядка относительно $\varphi$:
  \ba \label {f32}
  \lim_{\|\varphi\|\to0}\frac{\|D(\varphi,\mu)\|}{\|\varphi\|^k}=0.
  \ea
  
  При выполненных условиях уравнение
  \ba \label {f33}
  \varphi=A_1(\varphi,\mu)
  \ea
  имеет нулевое решение $\theta$ при всех значениях $\mu$.
  
  Интерес представляют случаи, когда наряду с тривиальным решением, при некотором значении $\mu$, близком к критическому\footnote{Такое значение параметра называется бифуркационным.} значению $\mu^*$, уравнение (\ref{f33}) имеет ещё одно (малое) решение, которое и будет определять форму волны.
  
  Развитая Красносельским теория бифуркаций (историю вопроса см. в [Богатов, Мухин, 2015]) включала в себя обоснование линеаризации для уравнения (\ref{f28}). Соответствующая теорема может быть сформулирована так [Красносельский, 1956, с.457]:
  
  \textbf{Теорема 1.} \textit{Пусть оператор $A_1(\varphi,\mu)$ представляется в виде (\ref{f30}) и выполнены условия (\ref{f31})-(\ref{f32}). Тогда каждая точка бифуркации этого оператора является характеристическим значением линейного оператора B. Каждое характеристическое значение нечётной кратности оператора B является точкой бифуркации оператора $A_1$.}
  
  Поскольку производная Фреше оператора $A_1$ в точке $\theta$ имеет вид 
  $$
  B\varphi=\int_0^{2\pi}K(s,t)\varphi(t)dt,
  $$
 [Красносельский, 1956a, c. 204], все его характеристические значения $\mu_n$ нечётнократны. Следовательно, все они будут точками бифуркации для операторов Некрасова (\ref{f28}).
  
  Отметим, что подход Красносельского позволял ответить на вопрос о точках бифуркации уравнения Некрасова \textit{без построения уравнения разветвления или исследования сходимости рядов} Назарова.
  
  Доказанные Красносельским теоремы о структуре множества решений уравнения (\ref{f33}) для вполне непрерывных операторов $A_1$ общего вида и его линеаризованных аналогов были основаны на свойствах вращения вполне непрерывных векторных полей [Красносельский, 1956a, Гл. IV]. Эти теоремы помогли обнаружить качественно новые свойства \textit{глобальных} решений уравнения Некрасова. В частности, было выяснено, что
  \begin{enumerate}
  	\item\textit{Решения $\varphi_\mu$ уравнения (\ref{f33}), отвечающие значениям параметра $\mu$, близким к $\mu_n$, образуют в окрестности $\theta$ непрерывную ветвь\footnote{Собственные функции некоторого оператора $A$ образуют в области $G\subset E$ \textit{непрерывную ветвь}, проходящую через точку $\varphi_0$, если граница каждого ограниченного открытого множества $U\subset G$, содержащего $\varphi_0$, имеет непустое пересечение c множеством собственных векторов $A$ [Красносельский, 1956, c. 190].}, причём $\mu\to\mu_n$, когда $\|\varphi_\mu\|\to 0$} [Красносельский, 1956a, c. 203].
   \item \textit{В каждой непрерывной ветви, существование которой обозначено в предыдущем пункте, фиксированным значениям $\mu$, близким к $\mu_n$, соответствует единственное ненулевое решение уравнения (\ref{f33}), непрерывно зависящее от $\mu$.} [Красносельский, 1956, c. 458].
 \end{enumerate}
   
\smallskip

{\bf 7. Переход к теории волн немалой амплитуды.}

 Топологические методы теории бифуркаций, разработанные Красносельским и основанные на степени отображения, позволили отказаться от предположения о малости амплитуды волны в теории Некрасова. Кроме того, поскольку физически осмысленные решения уравнений, подобных (\ref{f1}), (\ref{f12}) являются неотрицательными, для дальнейшего анализа уравнений вида (\ref{f33}) с оператором (\ref{f28}) имело смысл задействовать теорию конусов, основанную М.Г. Крейном и М.А. Рутманом [Крейн, Рутман, 1948] и развитую Красносельским (историю вопроса см. в [Богатов, 2020a]). Несмотря на то, что Красносельский хорошо ориентировался в задачах нелинейной механики (см., например, [Бахтин, Красносельский, 1955]), для дальнейшего качественного анализа задач о волнах установившегося вида лучше подходили специалисты, имеющие базовое образование в области механики и гидродинамики.
 
 Так сложилось, что в 1950-х гг. одна из школ нелинейной механики СССР была создана в Ростовском-на-Дону государственном университете (РГУ). Её организатором был выдающийся механик, впоследствии академик, Иосиф Израилевич Ворович. Он сам получил фундаментальную математическую подготовку (в том числе в области нелинейного анализа) на мехмате МГУ и хорошо понимал необходимость аналогичной подготовки для своих учеников, занимающихся вопросами нелинейной механики сплошной среды. По  словам Воровича [Ворович, Хапланов, 1963, с.213]
 
 \textit{«Коллектив механиков\footnote{Имеется в виду коллектив кафедры теоретической механики РГУ.} ощутил, что углубленная математическая разработка задач нелинейной механики возможна только на базе современных методов, в первую очередь, на базе нелинейного функционального анализа.»}
 
 В связи с этим в 1953 г. на кафедре теоретической механики РГУ был организован семинар по функциональному анализу, который имел тесные связи с одноимённым Воронежским семинаром. Это подтверждается выступлениями ростовчан на семинаре Красносельского [Красносельский, Крейн, Мышкис, 1957, с.249], а также выходом совместного сборника статей по функциональному анализу Воронежского госуниверситета и РГУ. Как отмечали Ворович и Хапланов [Ворович, Хапланов, 1963, с.213],\\ \textit{«Работа семинара сблизила ростовчан с воронежскими математиками, прежде всего с М.А. Красносельским, С.Г. Крейном и В.И. Соболевым.»}
 
 В результате один из учеников Воровича – Юрий Петрович Красовский – смог применить методы нелинейного функционального анализа к решению задач теории волн [Ворович, Хапланов, 1963, с. 226]. Красовский предполагал, что установившаяся жидкость бесконечной глубины движется без образования вихрей с постоянной скоростью, направленной по горизонтали, а профиль свободной границы является неподвижной периодической кривой. Отталкиваясь от постановки Леви-Чивита (\ref{f19}) и считая волны симметричными, Красовский пришёл к задаче отыскания периодической функции $\Phi(\varphi)\not\equiv \theta$ 
 \ba \label {f34}
 \Phi(\varphi)=\mu\int_0^\pi K_1(\varepsilon,\varphi)e^{3\Psi}\sin\Phi d\varepsilon,
 \ea
 где $\mu>0$ – число Фруда, а $K_1(\varepsilon,\varphi)$ – ядро, определяемое аналогично формуле (\ref{f29}) с заменой $\mu_n$ на величину $\frac{\pi n}{2}$. $\Psi$ – функция, сопряжённая с $\Phi$ и удовлетворяющая условию
 $$
 \int_0^{2\pi}\Psi(\varepsilon)d\varepsilon=0.
 $$
 
 Искомая функция $\Phi(\varphi)$  (как и в работах Некрасова) -- это угол, образованный касательной к профилю волны с горизонталью.
 
 Красовский ограничился случаем, когда волны имеют по одному гребню (точка $\varphi=\pi$) и одной впадине (точка $\varphi=0$). То есть 
   \begin{enumerate}
 	\item $\Phi(\varphi)\geq0$ при $\varphi\in[0,\pi]$;
 	\item $\Phi(0)=\Phi(\pi)=0$.
 \end{enumerate}

Это предположение позволило использовать теорию положительных операторов Красносельского [Красносельский, 1956a, Гл. V], краеугольным камнем которой является понятие \textit{конуса} $K$ банахова пространства $E$ – замкнутого выпуклого множества, содержащее вместе с каждой своей точкой $x$ луч, проходящий через неё и такого, что из того, что ${x,-x}\in K$ следует, что $x=\theta$. При помощи конуса $K$ в пространстве $E$ вводят полуупорядоченность, полагая $x\preceq y$, если $y-x\in K$(при этом знак $\preceq$ обладает всеми свойствами обычного знака $\leq$), а также определяют \textit{положительные операторы}, как операторы, оставляющие конус $K$ инвариантным. Большую пользу в контексте теории конусов приносит рассмотрение \textit{монотонных} операторов A, обладающими свойством $x\preceq y\Rightarrow Ax\preceq Ay$.

Вернёмся к уравнению (\ref{f34}) -- его можно записать в операторном виде:
 \ba \label {f35}
\Phi=\mu A\Phi,
\ea
где $A\Phi=\int_0^\pi K_1(\varepsilon,\varphi)e^{3\Psi}\sin\Phi d\varepsilon$.

Обозначим через $\dot C[0,\pi]$ – пространство непрерывных на $[0,\pi]$ функций, обращающихся в ноль на границе отрезка $[0,\pi]$; $K$ – конус неотрицательных в $\dot C[0,\pi]$ функций; $B_r$ – шар пространства $\dot C[0,\pi]$ радиуса $r<\frac\pi 6$. Тогда, как показал Красовский [Красовский, 1960, c. 1238], оператор $A$ из (\ref{f35}) будет вполне непрерывным в $B_r$ и положительным на множестве $K_r=B_r\cap K$. При этом оператор $A$ не является монотонным из-за наличия множителя $sin \, \Phi$. Тем не менее, поскольку он имеет монотонную миноранту (монотонный оператор $B$, для которого при всех $x\in K_r$ выполнено неравенство $Bx\preceq Ax$), к анализу уравнения (\ref{f35}) применимы теоремы о собственных функциях положительных операторов, доказанные Красносельским [Красносельский, 1956a, Гл. V, § 2-3]:

\textbf{Теорема 2.} \textit{Положительные собственные функции оператора A образуют в конусе K непрерывную ветвь длиной r}.

\textbf{Теорема 3.} \textit{Позитивный\footnote{\textit{Позитивный спектр} положительного оператора A – совокупность тех его собственных чисел, которым соответствуют положительные собственные функции.} спектр оператора A целиком лежит в некотором интервале $0<a<\mu <b$, причём a и b не зависят от $\|\Phi\|$.}

Из этих теорем следует, что в жидкости бесконечно большой глубины [Красовский, 1960, c. 1238]
 \begin{enumerate}
		\item при ограниченной скорости \textit{нет волн произвольно большой амплитуды};
		\item существуют периодические симметричные волны, у которых максимум угла наклона к профилю волны принимает любое значение из интервала $(0; \frac \pi 6)$.
\end{enumerate}

Похожие результаты были получены Красовским и при изучении волновых движений жидкости конечной глубины [Красовский, 1961].
   
\smallskip

{\bf Заключение.}

Появлению уравнения Некрасова предшествовали исследования двух ярких представителей школ механики: итальянской (Т. Леви-Чивиты, который стал изучать задачи об обтекании жидкостью препятствий методами теории функций комплексного переменного) и французской (А. Вилля, который получил нелинейное интегро-дифференциальное уравнение для формы обтекания жидкости с использованием подхода Леви-Чивиты). Наличие указанных результатов дало Некрасову возможность двигаться дальше и применить принятый в московской школе механики во главе с Жуковским подход, основанный на использовании функции Жуковского. Вкупе со взаимодействием с сильной математической школой МГУ во главе с Н.Н. Лузиным это помогло Некрасову перейти к более простому (интегральному) варианту уравнения, описывающему течение жидкости вокруг препятствия в задаче Вилля - Леви-Чивиты. Кроме того, это позволило распространить данный подход на задачи теории волн.

В силу ряда обстоятельств Некрасов, как и многие его коллеги (в том числе математики Н.Н. Лузин, А. Я. Хинчин и др.) был вынужден уехать из Москвы в 1920-м г. в г. Иваново в поисках лучшего жизнеобеспечения [Люстерник, 1967, с. 160]. До 1922 г. Некрасов был штатным профессором Иваново-Вознесенского политехнического института, организованного в 1918 г. [Секерж-Зенькович, 1960, c. 153]. Научная библиотека этого института, естественно, не шла ни в какое сравнение со столичными научными библиотеками, что затрудняло доступ профессуры к трудам своих зарубежных и даже петербургских коллег. Этим, по-видимому можно объяснить тот факт, что Некрасов не использовал  метод Ляпунова-Шмидта для решения своего уравнения, а создавал свой метод самостоятельно, «с нуля». С другой стороны, данное обстоятельство сослужило хорошую службу Некрасову, побуждая его к созданию теории нелинейных интегральных уравнений, ставшую хорошим дополнением к теории Ляпунова-Шмидта. Как выяснилось позднее,  теория Некрасова содержала в себе потенциал для её распространения на нелинейные уравнения общего вида в абстрактных пространствах. 

Поскольку методы, используемые Некрасовым – разложение в ряд по степеням малого параметра и мажорирование (как, впрочем, и метод Ляпунова-Шмидта) – имели локальный характер, ответить на вопросы о глубинных свойствах решений уравнения Некрасова они не могли. Здесь на помощь пришли топологические методы нелинейного анализа, развитые в работах М. А. Красносельского середины 1950-х гг. и основанные на теории степени отображения в банаховом пространстве. 

Это позволило выяснить 
\begin{itemize}
	\item сколько малых ненулевых решений имеет уравнение Некрасова при значениях параметра, близких к бифуркационному значению $\mu_n$;
	\item когда ненулевые решения образуют непрерывные ветви и сколько таких ветвей существует при $\mu\longrightarrow\mu_n$;
	\item как зависят от $\mu$ ненулевые решения уравнения Некрасова;
	\item при каких условиях на уравнение спектр оператора Некрасова будет сплошным.
\end{itemize}

Кроме того, применение Ю.П. Красовским теории конусов, разработанную Красносельским для решения задач нелинейного функционального анализа во второй половине 1950-х гг., позволило отказаться от предположения о малости амплитуды волны в уравнении Некрасова и провести его качественный анализ в новых условиях (1960). Немаловажную роль здесь сыграло тесное сотрудничество двух молодых научных школ СССР: воронежской школы функционального анализа (руководитель - М.А. Красносельский) и ростовской школы нелинейной механики (руководитель - И.И. Ворович).

Дальнейшее использование теории конусов применительно к глобальному анализу уравнения Некрасова было осуществлено в 1970-е гг. силами западных учёных [Kuznetsov, 2021, § 3.1].

Поводя итог, можно утверждать, что уравнение Некрасова достойно того, чтобы быть поставленным в один ряд с известными нелинейными уравнениями, таким как, например, уравнение Лиувилля, изучение которых даёт всё новые результаты и находит новое применение [Богатов, 2020b]. Ожидать ослабления интереса к уравнению Некрасова в ближайшем будущем не приходится.

\begin{center}
{\bf  Список литературы}
\end{center}

\begin{center}
{\bf References}
\end{center}

\end{document}